\documentclass[leqno,12pt]{amsart}
\usepackage{amssymb} 
\theoremstyle{plain} 

\theoremstyle{definition} 

\begin{document}

\title[Garnier system in two variables]{Remark on the Garnier system in two variables \\}
\author{Yusuke Sasano }

\renewcommand{\thefootnote}{\fnsymbol{footnote}}
\footnote[0]{2000\textit{ Mathematics Subjet Classification}.
34M55; 34M45; 58F05; 32S65.}

\maketitle

\begin{abstract}
We remark on the Garnier system in two variables.
\end{abstract}

\section{Summary}
In this note, we consider the following questions;

Why do we need Okamoto-Kimura's algebraic transformation of degree 2 for the Garnier system in two variables?

Here, the Garnier system in two variables is equivalent to the rational Hamiltonian system given by (see \cite{oka,10})
\begin{align}\label{1}
\begin{split}
dq_1&=\frac{\partial H_1}{\partial p_1}dt+\frac{\partial H_2}{\partial p_1}ds, \quad dp_1=-\frac{\partial H_1}{\partial q_1}dt-\frac{\partial H_2}{\partial q_1}ds,\\
dq_2&=\frac{\partial H_1}{\partial p_2}dt+\frac{\partial H_2}{\partial p_2}ds, \quad dp_2=-\frac{\partial H_1}{\partial q_2}dt-\frac{\partial H_2}{\partial q_2}ds,\\
H_1 &=-\frac{q_1(q_1-1)(q_1-t)(q_1-s)(q_2-t)}{(q_1-q_2)t(t-1)(t-s)}\{p_1^2+\frac{\kappa}{q_1(q_1-1)}\\
&-\left(\frac{\theta_1-1}{q_1-t}+\frac{\theta_2}{q_1-s}+\frac{\kappa_0}{q_1}+\frac{\kappa_1}{q_1-1} \right)p_1\}\\
&-\frac{q_2(q_2-1)(q_2-t)(q_2-s)(q_1-t)}{(q_2-q_1)t(t-1)(t-s)}\{p_2^2+\frac{\kappa}{q_2(q_2-1)}\\
&-\left(\frac{\theta_1-1}{q_2-t}+\frac{\theta_2}{q_2-s}+\frac{\kappa_0}{q_2}+\frac{\kappa_1}{q_2-1} \right)p_2\},\\
H_2&=\pi(H_1),
\end{split}
\end{align}
where the transformation $\pi$ is explicitly given by
\begin{align}
\begin{split}
\pi:&(q_1,p_1,q_2,p_2,t,s;\kappa_0,\kappa_1,\theta_1,\theta_2,\kappa) \rightarrow(q_2,p_2,q_1,p_1,s,t;\kappa_0,\kappa_1,\theta_2,\theta_1,\kappa).
\end{split}
\end{align}
Here, $q_1,p_1,q_2$ and $p_2$ are canonical variables and $\kappa_0,\kappa_1,\theta_1,\theta_2$ and $\kappa$ are constant parameters satisfying the relation
\begin{equation}
\kappa=\frac{1}{4}[(\kappa_0+\kappa_1+\theta_1+\theta_2-1)^2-\kappa_{\infty}^2].
\end{equation}

For the system \eqref{1}, we calculate its symmetry. We show that each B{\"a}cklund transformation is a coupled B{\"a}cklund transformation of the Painlev\'e VI system.

We see that the system \eqref{1} is invariant under the following transformations defined as follows\rm{:\rm} with the notation $(*)=(q_1,p_1,q_2,p_2,t,s;\kappa_0,\kappa_1,\kappa_{\infty},\theta_1,\theta_2),$
\begin{align}
\begin{split}
        s_0: (*) \rightarrow &(\frac{1}{q_1},-\left(q_1p_1-\frac{\kappa_0+\kappa_1-\kappa_{\infty}+\theta_1+\theta_2-1}{2} \right)q_1,\\
        &\frac{1}{q_2},-\left(q_2p_2-\frac{\kappa_0+\kappa_1-\kappa_{\infty}+\theta_1+\theta_2-1}{2} \right)q_2,\frac{1}{t},\frac{1}{s};\\
        &-\kappa_{\infty},\kappa_1,\kappa_0,\theta_1,\theta_2),\\
        s_1: (*) \rightarrow &\left(q_1,p_1-\frac{\kappa_0}{q_1},q_2,p_2-\frac{\kappa_0}{q_2},t,s;-\kappa_0,\kappa_1,\kappa_{\infty},\theta_1,\theta_2 \right),\\
        s_2: (*) \rightarrow &\left(q_1,p_1-\frac{\kappa_1}{q_1-1},q_2,p_2-\frac{\kappa_1}{q_2-1},t,s;\kappa_0,-\kappa_1,\kappa_{\infty},\theta_1,\theta_2 \right),\\
        s_3: (*) \rightarrow &\left(q_1,p_1,q_2,p_2,t,s;\kappa_0,\kappa_1,-\kappa_{\infty},\theta_1,\theta_2 \right),\\
        s_4: (*) \rightarrow &\left(q_1,p_1-\frac{\theta_1}{q_1-t},q_2,p_2-\frac{\theta_1}{q_2-t},t,s;\kappa_0,\kappa_1,\kappa_{\infty},-\theta_1,\theta_2 \right),\\
        s_5: (*) \rightarrow &\left(q_1,p_1-\frac{\theta_2}{q_1-s},q_2,p_2-\frac{\theta_2}{q_2-s},t,s;\kappa_0,\kappa_1,\kappa_{\infty},\theta_1,-\theta_2 \right),\\
        {\sigma}_1: (*) \rightarrow &(1-q_1,-p_1,1-q_2,-p_2,1-t,1-s;\kappa_1,\kappa_0,\kappa_{\infty},\theta_1,\theta_2),\\
        {\sigma}_2: (*) \rightarrow &(q_2,p_2,q_1,p_1,s,t;\kappa_0,\kappa_1,\kappa_{\infty},\theta_2,\theta_1),\\
        {\sigma}_3: (*) \rightarrow &\left(\frac{s-q_1}{s-1},-(s-1)p_1,\frac{s-q_2}{s-1},-(s-1)p_2,\frac{s-t}{s-1},\frac{s}{s-1};\theta_2,\kappa_1,\kappa_{\infty},\theta_1,\kappa_0 \right),\\
        {\sigma}_4: (*) \rightarrow &(\frac{1}{q_1},-\left(q_1p_1-\frac{\kappa_0+\kappa_1+\kappa_{\infty}+\theta_1+\theta_2-1}{2} \right)q_1,\\
        &\frac{1}{q_2},-\left(q_2p_2-\frac{\kappa_0+\kappa_1+\kappa_{\infty}+\theta_1+\theta_2-1}{2} \right)q_2,\frac{1}{t},\frac{1}{s};\\
        &\kappa_{\infty},\kappa_1,\kappa_0,\theta_1,\theta_2).
        \end{split}
        \end{align}
The group $<{\sigma}_1,{\sigma}_2,{\sigma}_3,{\sigma}_4>$ is isomorphic to symmetric group of degree five.

By resolving an accessible singularity of the system \eqref{1}, we transform the system \eqref{1} into a polynomial Hamiltonian system.

We see that the birational and symplectic transformation $\varphi_1$\rm{:\rm}
\begin{equation}\label{symp}
  \left\{
  \begin{aligned}
   Q_1=&\frac{1}{q_1-q_2},\\
   P_1=&-\left((q_1-q_2)p_1-\frac{\kappa_0+\kappa_1-\kappa_{\infty}+\theta_1+\theta_2-1}{2}\right)(q_1-q_2),\\
   Q_2=&q_2,\\
   P_2=&p_2+p_1
   \end{aligned}
  \right. 
\end{equation}
takes the system \eqref{1} into a polynomial Hamiltonian system.

We remark that for the polynomial Hamiltonian system obtained by \eqref{symp}, this system becomes again a polynomial Hamiltonian system in each coordinate $r_i \ (i=0,1,\ldots,6)${\rm : \rm}
\begin{align*}
&r_0:x_0=q_1, \ y_0=p_1, \ z_0=-(q_2p_2-\kappa_0)p_2, \ w_0=\frac{1}{p_2}, \\
&r_1:x_1=q_1, \ y_1=p_1, \ z_1=-((q_2-1)p_2-\kappa_1)p_2, \ w_1=\frac{1}{p_2}, \\
&r_2:x_2=(q_1q_2+1)q_2, \ y_2=\frac{p_1}{q_2^2}, \ z_2=\frac{1}{q_2}, \ w_2=-\left(q_2p_2-2\left(q_1q_2+\frac{1}{2} \right)\frac{p_1}{q_2} \right)q_2, \\
&r_3:x_3=-(q_1p_1+\kappa_{\infty})p_1, \ y_3=\frac{1}{p_1}, \ z_3=q_2, \ w_3=p_2,\\
&r_4:x_4=q_1, \ y_4=p_1, \ z_4=-((q_2-t)p_2-\theta_1)p_2, \ w_4=\frac{1}{p_2},\\
&r_5:x_5=q_1, \ y_5=p_1, \ z_5=-((q_2-s)p_2-\theta_2)p_2, \ w_5=\frac{1}{p_2},\\
&r_6:x_6=q_1, \ y_6=p_1+\frac{p_2}{q_1^2}, \ z_6=q_2+\frac{1}{q_1}, \ w_6=p_2.
\end{align*}
Here, for notational convenience, we have renamed $Q_i,P_i$ to $q_i,p_i$ (which are not the same as the previous $q_i,p_i$).

We see that the system obtained by \eqref{symp} is invariant under the following transformations defined as follows\rm{:\rm} with the notation $(*)=(q_1,p_1,q_2,p_2,t,s;\kappa_0,\kappa_1,\kappa_{\infty},\theta_1,\theta_2);$
\begin{align}
\begin{split}
    \varphi_1 \circ  s_1 \circ s_0 \circ {\varphi_1}^{-1}: (*) \rightarrow &(-(q_1q_2+1)q_2,-\frac{p_1}{q_2^2},\frac{1}{q_2},-\left(q_2p_2-2\left(q_1q_2+\frac{1}{2} \right)\frac{p_1}{q_2} \right)q_2,\frac{1}{t},\frac{1}{s};\\
        &-\kappa_{\infty},\kappa_1,-\kappa_0,\theta_1,\theta_2),\\
        \varphi_1 \circ s_1 \circ {\varphi_1}^{-1}: (*) \rightarrow &\left(q_1,p_1-\frac{\kappa_0 q_2}{q_1 q_2+1},q_2,p_2-\frac{\kappa_0  (2 q_1 q_2+1)}{q_2(q_1 q_2+1)},t,s;-\kappa_0,\kappa_1,\kappa_{\infty},\theta_1,\theta_2 \right),\\
        \varphi_1 \circ {\sigma}_2 \circ {\varphi_1}^{-1}: (*) \rightarrow &(-q_1,-\left(p_1+\frac{p_2}{q_1^2} \right),q_2+\frac{1}{q_1},p_2,s,t;\kappa_0,\kappa_1,\kappa_{\infty},\theta_2,\theta_1).
        \end{split}
        \end{align}
On the other hand, the system obtained by \eqref{symp} is not invariant under the following transformation associated with holomorphy condition $r_0$:
\begin{align}
\begin{split}
        S_0: (*) \rightarrow &\left(q_1,p_1,q_2,p_2-\frac{\kappa_0}{q_2},t,s;-\kappa_0,\kappa_1,\kappa_{\infty},\theta_1,\theta_2 \right).
        \end{split}
        \end{align}
In this vein, we will guess that in \cite{Oka} Professors H. Kimura and K. Okamoto considered an algebraic transformation of degree 2 for the Garnier system in two variables.

\end{document}